\documentclass[12pt]{article}

\usepackage{authblk}
\usepackage{pgfplots}
\usepackage{tikz}

\usetikzlibrary{automata}
\usetikzlibrary{positioning}
\usepackage{amsmath, amssymb} 
\usepackage{amsthm}
\usepackage{graphicx}
\usepackage[shortlabels]{enumitem}
\usepackage{caption}

\usepackage{txfonts,geometry,hyperref}

\newtheorem{conj}{Conjecture}[section]

\newtheorem{theorem}{Theorem}[section]

\newtheorem{cor}[theorem]{Corollary}

\DeclareMathOperator{\Spec}{Spec}
\DeclareMathOperator{\L-Spec}{L-spec}
\DeclareMathOperator{\Q-Spec}{Q-spec}
\DeclareMathOperator{\CN-Spec}{CN-spec}

\begin{document}
\title{Characterization of  commuting graphs of finite groups having small genus }
\author{Shrabani Das$^1$, Deiborlang Nongsiang$^2$ and Rajat Kanti Nath$^3$\footnote{Corresponding author}}
\affil{$^{1, 3}$Department of Mathematical Sciences, Tezpur University, Napaam-784028, Sonitpur, Assam, India.

$^2$Department of Mathematics, North-Eastern Hill University, Shillong-793022, Meghalaya, India.

Emails: shrabanidas904@gmail.com (S. Das); ndeiborlang@yahoo.in (D. Nongsiang); rajatkantinath@yahoo.com (R. K. Nath)}
\date{ }
\maketitle

\begin{abstract}
In this paper we first show that among all double-toroidal and triple-toroidal finite graphs only $K_8 \sqcup 9K_1$, $K_8 \sqcup 5K_2$, $K_8 \sqcup 3K_4$, $K_8 \sqcup 9K_3$, $K_8\sqcup 9(K_1 \vee 3K_2)$, $3K_6$ and $3K_6 \sqcup 4K_4 \sqcup 6K_2$ can be realized as commuting graphs of finite groups. As consequences of our results we also show that  for any finite non-abelian  group $G$ if  the commuting graph of $G$ (denoted by $\Gamma_c(G)$) is double-toroidal or  triple-toroidal then $\Gamma_c(G)$ and its complement satisfy Hansen-Vuki{\v{c}}evi{\'c} Conjecture and E-LE conjecture. 	
In the process we  find a non-complete graph, namely the non-commuting graph of the group $(\mathbb{Z}_3 \times \mathbb{Z}_3) \rtimes Q_8$, that is hyperenergetic. This gives a new counter example to a conjecture of Gutman regarding hyperenergetic graphs.
\end{abstract}
\section{Introduction}
     Finite groups are being characterized through various graphs defined on it for a long time now. A survey on graphs defined on groups can be found in \cite{cameron2021graphs}. One such graph defined on groups is the commuting graph. The commuting graph of  a finite group was originated from the works of Brauer and Fowler in \cite{brauer1955groups}. Let $G$ be a finite non-abelian group with center $Z(G)$. The commuting graph of $G$ is a simple undirected graph whose vertex set is $G \setminus Z(G)$ and two vertices $x$ and $y$ are adjacent if  $xy=yx$. It is denoted by $\Gamma_{c}(G) $. The complement of this graph is the non-commuting graph of $G$, denoted by $ \Gamma_{nc}(G)$. The study of non-commuting graph of a finite non-abelian group gets popularity because of a question posed by Erd{\"o}s in the year 1975 which was answered by Neumann in 1976 \cite{neumann1976problem}.

      The genus of a graph $\Gamma$ is the smallest non-negative integer $n$ such that the graph can be embedded on the surface obtained by attaching $n$ handles to a sphere. It is denoted by $\gamma(\Gamma)$. The graphs which have genus zero are called planar graphs, those which have genus one are called toroidal graphs, those which have genus two are called double-toroidal graphs and those which have genus three are called triple-toroidal graphs. Classification of finite non-abelian groups whose commuting graphs are planar or toroidal can be found in \cite{AFK2015} and \cite{DN-2016} (also see \cite[Theorem 3.3]{B1}). Recently, finite non-abelian groups such that their commuting graphs are double-toroidal or triple-toroidal are classified in  \cite{DNongsiang}.  In this paper, we consider finite non-abelian groups whose commuting graphs are double or triple-toroidal and realize their commuting graphs. As such we 
      show that among all double-toroidal and triple-toroidal finite graphs only $K_8 \sqcup 9K_1$, $K_8 \sqcup 5K_2$, $K_8 \sqcup 3K_4$, $K_8 \sqcup 9K_3$, $K_8\sqcup 9(K_1 \vee 3K_2)$, $3K_6$ and $3K_6 \sqcup 4K_4 \sqcup 6K_2$ can be realized as commuting graphs of finite groups. We also compute  first and second Zagreb indices of $\Gamma_c(G)$ and $\Gamma_{nc}(G)$  and show that   they satisfy Hansen-Vuki{\v{c}}evi{\'c} conjecture if $\Gamma_c(G)$ is double-toroidal or triple-toroidal. Further, we show that these graphs also satisfy E-LE conjecture.


        Let $\Gamma$ be a simple undirected graph with vertex set $v(\Gamma)$ and edge set $e(\Gamma)$. The first and second Zagreb indices of $\Gamma$, denoted by $M_{1}(\Gamma)$ and $M_{2}(\Gamma)$ respectively, are defined as 
\[
M_{1}(\Gamma) = \sum\limits_{v \in v(\Gamma)} \deg(v)^{2}  \text{ and }  M_{2}(\Gamma) = \sum\limits_{uv \in e(\Gamma)} \deg(u)\deg(v),
\]
where $ \deg(v) $ is the number of edges incident on $ v $
(called degree of $v$). Zagreb indices of graphs were introduced by Gutman and Trinajsti{\'c} \cite{Gut-Trin-72} in 1972 to examine the dependence of total $\pi$-electron energy on molecular structure. As noted in \cite{Z-index-30y-2003}, Zagreb indices are also used in studying molecular complexity, chirality,    ZE-isomerism and heterosystems etc. Later on, general mathematical properties of these indices  are also studied by many mathematicians. A survey on mathematical properties of Zagreb indices can be found in \cite{Gut-Das-2004}.
 Comparing first and second Zagreb indices, Hansen and Vuki{\v{c}}evi{\'c} \cite{hansen2007comparing} posed the following   conjecture in  2007.
\begin{conj}\label{Conj}
(Hansen-Vuki{\v{c}}evi{\'c} Conjecture) For any simple finite graph $\Gamma$, 
\begin{equation}\label{Conj-eq}
 \dfrac{M_{2}(\Gamma)}{\vert e(\Gamma) \vert} \geq \dfrac{M_{1}(\Gamma)}{\vert v(\Gamma) \vert} .
\end{equation}
\end{conj}
It was shown in \cite{hansen2007comparing} that the conjecture is not true  if $\Gamma = K_{1, 5} \sqcup K_3$. However, Hansen and Vuki{\v{c}}evi{\'c} \cite{hansen2007comparing}  showed that  Conjecture \ref{Conj}  holds for chemical graphs. In \cite{vukicevic2007comparing}, it was shown that the conjecture holds for trees with equality in \eqref{Conj-eq} when $\Gamma$ is a star graph. In \cite{liu2008conjecture}, it was shown that the conjecture holds for connected unicyclic graphs with equality when the graph is a cycle. However, the search of graphs validating or invalidating Conjecture \ref{Conj} is not completed yet. Recently, Das et al. \cite{SD-AS-RKN-2023} have obtained various finite non-abelian groups such that their commuting graphs satisfy Hansen-Vuki{\v{c}}evi{\'c} Conjecture. It was also shown that $\Gamma_c(G)$  satisfies Hansen-Vuki{\v{c}}evi{\'c} Conjecture if $\Gamma_c(G)$ is planar or toroidal.

Let $A(\Gamma)$ and $D(\Gamma)$ denote the adjacency matrix and degree matrix of $\Gamma$ respectively. The set of eigenvalues of $A(\Gamma)$ along with their multiplicities is called the spectrum of $\Gamma$. The Laplacian matrix and signless Laplacian matrix of $\Gamma$ are given by $L(\Gamma):=D(\Gamma)-A(\Gamma)$ and $Q(\Gamma):=D(\Gamma)+A(\Gamma)$ respectively. The Laplacian spectrum and signless Laplacian spectrum are the set of eigenvalues of $L(\Gamma)$ and $Q(\Gamma)$ along with their multiplicities respectively. Let $v(\Gamma):=\{v_i:i=1,2, \ldots, n\}$. The common neighbourhood of two distinct vertices $v_i$ and $v_j$, denoted by $C(v_i, v_j)$, is the set of all vertices other than $v_i$ and $v_j$ which are adjacent to both $v_i$ and $v_j$. The common neighbourhood matrix of $\Gamma$, denoted by $CN(\Gamma)$, is defined as
\[
(CN(\Gamma))_{i, j}= \begin{cases}
	|C(v_i, v_j)|, & \text{if } i \neq j \\
	0, & \text{if } i=j.
\end{cases}
\]
The common neighbourhod spectrum of $\Gamma$ is the set of all eigenvalues of $CN(\Gamma)$ along with their multiplicities. We write $\Spec(\Gamma)$, $\L-Spec(\Gamma)$, $\Q-Spec(\Gamma)$ and $\CN-Spec(\Gamma)$ to denote the spectrum, Laplacian spectrum, signless Laplacian spectrum and common neighbourhood spectrum of $\Gamma$ respectively. 

The energy, $E(\Gamma)$ and common neighbourhood energy, $E_{CN}(\Gamma)$ of $\Gamma$ are the sum of the absolute values of the elements of 
$\Spec(\Gamma)$ and $\CN-Spec(\Gamma)$ respectively. The Laplacian energy, $LE(\Gamma)$ and signless Laplacian energy, $LE^{+}(\Gamma)$ of $\Gamma$ are defined as
\[
LE(\Gamma)=\sum_{\lambda \in \L-Spec(\Gamma)} \left| \lambda-\frac{2m}{n} \right| \quad \text{and} \quad LE^{+}(\Gamma)= \sum_{\mu \in \Q-Spec(\Gamma)} \left| \mu-\frac{2m}{n} \right|,
\]
where $m = |e(\Gamma)|$. It is well known that $E(K_n)=LE(K_n)=LE^+(K_n)=2(n-1)$ and $E_{CN}(K_n)=2(n-1)(n-2)$.
A graph $\Gamma$ with $|v(\Gamma)| =n$  is called hyperenergetic if $E(\Gamma) > E(K_n)$. It is called hypoenergetic if $E(\Gamma) < n$. Similarly, $\Gamma$ is called L-hyperenergetic if $LE(\Gamma) > LE(K_n)$, Q-hyperenergetic if $LE^+(\Gamma) > LE^+(K_n)$ and CN-hyperenergetic if $E_{CN}(\Gamma) > E_{CN}(K_n)$.

 Gutman et al. \cite{E-LE-Gutman}  conjectured that $E(\Gamma) \leq LE(\Gamma)$ which is known as E-LE conjecture. Gutman \cite{Gutman-78} also conjectured that ``$\mathcal{G}$ is not hyperenergetic if $\mathcal{G} \ncong K_{|v(\mathcal{G})|}$". Note that both the conjectures were disproved. However, it is still unknown whether  the commuting or non-commuting graphs of finite groups satisfy E-LE conjecture. In this paper, we show that $\Gamma_{c}(G)$ and $\Gamma_{nc}(G)$ satisfy E-LE conjecture if $\Gamma_{c}(G)$ is double-toroidal or triple-toroidal. Further, we  find a non-complete graph, namely the non-commuting graph of the group $(\mathbb{Z}_3 \times \mathbb{Z}_3) \rtimes Q_8$, that is hyperenergetic. This gives a new counter example to the above mentioned conjecture of Gutman. We shall also determine whether $\Gamma_{c}(G)$ and $\Gamma_{nc}(G)$ are hypoenergetic, hyperenergetic, L-hyperenergetic, Q-hyperenergetic and CN-hyperenergetic if $\Gamma_{c}(G)$ is double-toroidal or triple-toroidal.

\section{Realization of commuting graph}
In this section, we determine all finite planar, toroidal, double-toroidal and triple-toroidal graphs that can be realized as commuting graphs of finite groups. Using  \cite[Theorem 2.2]{AFK2015}, \cite[Theorem 3.3]{B1} and commuting graphs of various finite non-abelian groups considered in  \cite{B1} we have the following theorem.
\begin{theorem}
	\begin{enumerate}[\rm (a)]
		\item Among all the planar finite graphs only $K_2 \sqcup 3K_1, 3K_2, K_4 \sqcup 5K_1, K_4 \sqcup 3K_2, 3K_4, K_3 \sqcup 4K_2, 5K_3 \sqcup 10K_2 \sqcup 6K_4, 3K_2 \sqcup 4K_4, K_4 \sqcup 5K_3$ and $7K_2 \sqcup D$ can be realized as commuting graphs of finite groups, where $D$ is the graph obtained from $4K_3$ after three vertex contractions as shown in Figure \ref{fig:fig1}.
		\item Among all the toroidal finite graphs only $K_6 \sqcup 7K_1, K_6 \sqcup 4K_2, K_6 \sqcup 3K_3, K_6 \sqcup 4K_4$ and $K_6 \sqcup 7K_2$ can be realized as commuting graphs of finite groups.
	\end{enumerate}
\end{theorem}
 \vspace{0.5cm} 
  \begin{minipage}[t]{0.5\textwidth}
     \begin{tikzpicture}
	\tikzstyle{vertex}=[circle,minimum size=0.1pt,fill,inner sep=2pt]
	\node[vertex](A) at (0,0){};
	\node[vertex](B) at (-1,-1){};
	\node[vertex](C) at (1,-1){};
	\node[vertex](D) at (3,-1){};
	\node[vertex](E) at (2,-2){};
	\node[vertex](F) at (1,-3){};
        \node[vertex](G) at (3,-3){};
        \node[vertex](H) at (4,0){};
        \node[vertex](I) at (5,-1){};
    \path (A) edge (B)
          (A) edge (C)
          (B) edge (C)
          (C) edge (D)
          (C) edge (E)
          (D) edge (E)
          (E) edge (F)
          (E) edge (G)
          (F) edge (G)
          (D) edge (H)
          (D) edge (I)
          (H) edge (I);
          \end{tikzpicture}
     
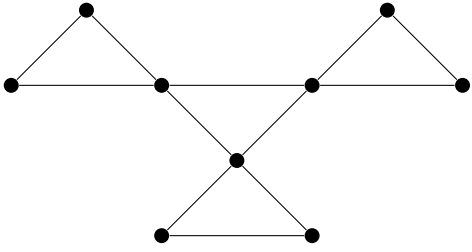
\captionof{figure}{Graph after three vertex\\ contractions in $4K_3$}
         \label{fig:fig1}
 \end{minipage}
\begin{minipage}[t]{0.5\textwidth}
    \begin{tikzpicture}
 \tikzstyle{vertex}=[circle,fill=black,inner sep=2pt]
 \node[vertex](A) at (1,0){};
 \node[vertex](F) at (3,0){};
 \node[vertex](B) at (0,-1){};

\node[vertex](G) at (4,-1){};
\node[vertex](C) at (2,-1) {};
 \node[vertex](D) at (1,-2) {};
 \node[vertex](E) at (3,-2){};

 \path 
 (A) edge (B)
 (B) edge (C)
 (A) edge (C)
 (C) edge (D)
 (D) edge (E)
 (C) edge (E)
 (C) edge (F)
 (G) edge (C)
 (F) edge (G);
\end{tikzpicture}	

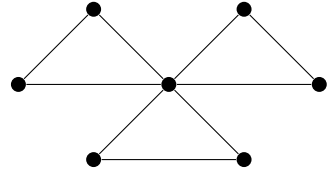
\captionof{figure}{$K_1 \vee 3K_2$~~~~~~~~~~~~~~~~~~~~~~~~~~~}
\label{fig:fig2}
\end{minipage}

%

 \vspace{.5cm} 

 The  following two results from \cite{DNongsiang} are useful in determining all finite double-toroidal and triple-toroidal graphs that can be realized as commuting graphs of finite groups.
\begin{theorem}\label{groups for double toroidal}
\cite{DNongsiang} Let $G$ be a finite non-abelian group. Then the commuting graph of $G$ is double-toroidal if and only if $G$ is isomorphic to one of the following groups:
    \begin{enumerate}[\rm (a)]
        \item $D_{18}, D_{20}, Q_{20}, S_3 \times \mathbb{Z}_{2} \times \mathbb{Z}_{2}, S_3 \times \mathbb{Z}_{4}$,
        \item $(\mathbb{Z}_3 \times \mathbb{Z}_3) \rtimes \mathbb{Z}_{2} \cong \langle x, y, z : x^3=y^3=z^2=[x, y]=1, x^z=x^{-1}, y^z=y^{-1} \rangle$,
        \item $\mathbb{Z}_3 \rtimes \mathbb{Z}_8 \cong \langle x, y : x^8=y^3=1, y^x=y^{-1} \rangle$,
        \item $(\mathbb{Z}_3 \rtimes \mathbb{Z}_4) \times \mathbb{Z}_2 \cong \langle x, y, z : x^4=y^3=z^2=1, xyx^{-1}=y^{-1}, xz=zx, yz=zy \rangle$, 
        \item $(\mathbb{Z}_3 \times \mathbb{Z}_3) \rtimes \mathbb{Z}_4 \cong \langle x, y : x^4=y^3=(yx^2)^2=[x^{-1}yx, y]=1 \rangle$, 
        \item $(\mathbb{Z}_3 \times \mathbb{Z}_3) \rtimes Q_8 \cong \langle x, y, z : x^4=y^4=z^3=1, y^x=y^{-1}, z^{y^2}=z^{-1}, z^{x^2}=z^{-1}, x^{-1}zx^{-1}=(zy)^2 \rangle$.
    \end{enumerate}
\end{theorem}
\begin{theorem}\label{groups for triple toroidal}
\cite{DNongsiang} Let $G$ be a finite non-abelian group. Then the commuting graph of $G$ is triple-toroidal if and only if $G$ is isomorphic to either
    \begin{enumerate}[\rm (a)]
        \item $GL(2, 3), D_8 \times \mathbb{Z}_3, Q_8 \times \mathbb{Z}_3$,
        \item $SL(2,3) \circ \mathbb{Z}_2\cong \langle x,y,z: y^3=z^4=1, x^2=z^2, y^x=y^{-1}, y^{-1}zy^{-1}z^{-1}y^{-1}z= xz^{-1}xy^{-1}zy=1 \rangle.$
    \end{enumerate}
\end{theorem}


Now we realize the structures of $\Gamma_c(G)$ if  $\Gamma_c(G)$ is double-toroidal or triple-toroidal.
\begin{theorem}\label{graphs when double toroidal}
    Let $G$ be a finite non-abelian group. If \, $\Gamma_c(G)$ is double-toroidal then $\Gamma_c(G)$ is isomorphic to $K_8 \sqcup 9K_1$, $K_8 \sqcup 5K_2$, $K_8 \sqcup 3K_4$, $K_8 \sqcup 9K_3$ or $K_8\sqcup 9(K_1 \vee 3K_2)$.
\end{theorem}
\begin{proof}

    From Theorem \ref{groups for double toroidal}, we have $\Gamma_c(G)$ is double-toroidal if and only if $G$ is isomorphic to either $D_{18}, D_{20}, Q_{20}, S_3 \times \mathbb{Z}_{2} \times \mathbb{Z}_{2}, S_3 \times \mathbb{Z}_{4}$, $(\mathbb{Z}_3 \times \mathbb{Z}_3) \rtimes \mathbb{Z}_{2}$, $\mathbb{Z}_3 \rtimes \mathbb{Z}_8$, $(\mathbb{Z}_3 \rtimes \mathbb{Z}_4) \times \mathbb{Z}_2$, $(\mathbb{Z}_3 \times \mathbb{Z}_3) \rtimes \mathbb{Z}_4$ or $(\mathbb{Z}_3 \times \mathbb{Z}_3) \rtimes Q_8$. 
    
    Let $G$ be any of the groups $D_{18}$ and $(\mathbb{Z}_3\times \mathbb{Z}_3)\rtimes\mathbb{Z}_2$. Then $G$ is an  AC-groups. The centralizers of the non-central elements of $G$ are of size 9 and 2. There is exactly one centralizer of size 9  and nine distinct centralizers of size 2. Thus $\Gamma_c(G)\cong K_8 \sqcup 9K_1$.

    Let $G$ be any of the groups $D_{20}$ and $Q_{20}$. Then $G$ is an AC-group, $|Z(G)|=2$ and it has one centralizer of size 10 and 5 distinct centralizers of size 4. Thus $\Gamma_c(G)\cong K_8 \sqcup 5K_2$.
    










    
    Let $G$ be any of the groups $S_3 \times \mathbb{Z}_{2} \times \mathbb{Z}_{2}, S_3 \times \mathbb{Z}_{4}, \mathbb{Z}_3 \rtimes \mathbb{Z}_8$ and $(\mathbb{Z}_3 \rtimes \mathbb{Z}_4) \times \mathbb{Z}_2$. Then $G$ is an AC-group, $|Z(G)|=4$ and it has three centralizer of size 8 and one centralizer of size 12. Thus $\Gamma_c(G)\cong K_8 \sqcup 3K_4$.

    If $G=(\mathbb{Z}_3 \times \mathbb{Z}_3) \rtimes \mathbb{Z}_4$, then $G$ is an AC-group, $|Z(G)|=1$ and it has one centralizer of size 9 and 9 centralizers of size 4. Thus $\Gamma_c(G)\cong K_8 \sqcup 9K_3$.

    Let $G=(\mathbb{Z}_3 \times \mathbb{Z}_3) \rtimes Q_8$. The group $G$ consist of one sylow 3-subgroup of order 9 and 9 sylow 2-subgroups of order 8. The sylow 2-subgroups of $G$ are isomorphic to $Q_8$ and the sylow 3-subgroup is isomorphic to $\mathbb{Z}_3 \times \mathbb{Z}_3$. The intersection of any two of these subgroups is trivial. Thus $G$ is exactly the union of these subgroups. Let $L$ be any of these subgroups and $x\in L, x\neq 1$. Then $C_{G}(x)\subseteq L$. Thus the commuting graph of $G$ consist of 10 components. One of the component is $\Gamma_c(G)[H]$, where $H\cup \{1\}$ is the sylow 3-subgroup of $G$. The other 9 components are $\Gamma_c(G)[K_i]$, where $K_i\cup \{1\},i=1,2,\dots,9,$ are the sylow 2-subgroups of $G$. It can be seen that $\Gamma_c(G)[H]\cong K_8$ and 
    $\Gamma_c(G)[K_i]\cong K_1 \vee 3K_2$ for $i=1,2,\dots,9$. Thus $\Gamma_c(G)\cong K_8\sqcup 9(K_1 \vee 3K_2)$.
\end{proof}

\begin{theorem}\label{graphs when triple toroidal}
    Let $G$ be a finite non-abelian group. If \, $\Gamma_c(G)$ is triple-toroidal, then $\Gamma_c(G)$ is isomorphic to  $3K_6$ or $3K_6 \sqcup 4K_4 \sqcup 6K_2$.
\end{theorem}
\begin{proof}
    From Theorem \ref{groups for triple toroidal}, we have $\Gamma_c(G)$ is triple-toroidal if and only if $G$ is isomorphic to  $GL(2, 3), D_8 \times \mathbb{Z}_3, Q_8 \times \mathbb{Z}_3$ or $C_2\circ S_4$.
    
    If $G = D_8 \times \mathbb{Z}_3$ or $Q_8 \times \mathbb{Z}_3$, then $G$ is an AC-group, $|Z(G)|=6$ and has three distinct centralizers of size 12. Therefore, $\Gamma_c(G)=K_6 \sqcup K_6 \sqcup K_6= 3K_6$.

   If $G = GL(2, 3)$ or $C_2\circ S_4$, then $G$ is an AC-group, $|Z(G)|=2$ and it has three centralizers of size 8, four centralizers of size 6 and six centralizers of size 4. Thus $\Gamma_c(G)=3K_6 \sqcup 4K_4 \sqcup 6K_2$.
\end{proof}

We conclude this section with the following corollary.
\begin{cor}
\begin{enumerate}[\rm (a)]
	\item Among all the double-toroidal finite graphs only $K_8 \sqcup 9K_1$, $K_8 \sqcup 5K_2$, $K_8 \sqcup 3K_4$, $K_8 \sqcup 9K_3$ or $K_8\sqcup 9(K_1 \vee 3K_2)$ can be realized as commuting graphs of finite groups.	
	\item Among all the triple-toroidal finite graphs only $3K_6$ and $3K_6 \sqcup 4K_4 \sqcup 6K_2$ can be realized as commuting graphs of finite groups.
\end{enumerate}
\end{cor}

\section{Some consequences}
In this section we show that  for any finite non-abelian  group $G$ if  $\Gamma_c(G)$ is double-toroidal or $\Gamma_c(G)$ is triple-toroidal then $\Gamma_c(G)$ and $\Gamma_{nc}(G)$ satisfy Hansen-Vuki{\v{c}}evi{\'c} Conjecture and E-LE conjecture. The following result is useful in our study.
\begin{theorem}\label{using commuting finding complement}
	(\cite{Das-Kex-14}, Page 575 and \cite{Das-Gut-04}, Lemma 3) For any graph  $\Gamma$ and its complement $\overline{\Gamma}$,
	\[
	M_{1}(\overline{\Gamma})=|v(\Gamma)|(|v(\Gamma)|-1)^{2}-4|e(\Gamma)|(|v(\Gamma)|-1)+ M_{1}(\Gamma) \quad  \text{ and }
	\]
	\[ 
	M_{2}(\overline{\Gamma})= \dfrac{|v(\Gamma)|(|v(\Gamma)|-1)^{3}}{2}+2|e(\Gamma)|^{2}-3|e(\Gamma)|(|v(\Gamma)|-1)^{2}+\left(|v(\Gamma)|-\dfrac{3}{2}\right)M_{1}(\Gamma)-M_{2}(\Gamma). 
	\]
\end{theorem}

In \cite{SD-AS-RKN-2023}, it was shown that $\Gamma_c(G)$  satisfies Hansen-Vuki{\v{c}}evi{\'c} Conjecture if $\Gamma_c(G)$ is planar or toroidal. 
The following theorem shows that if  $\Gamma_c(G)$ is double-toroidal then    $\Gamma_c(G)$ and $\Gamma_{nc}(G)$ satisfy Hansen-Vuki{\v{c}}evi{\'c} Conjecture.  
\begin{theorem}
    Let $G$ be a finite non-abelian group. If \, $\Gamma_c(G)$ is double-toroidal then $\frac{M_{2}(\Gamma(G))}{|e(\Gamma(G))|} \geq \frac{M_{1}(\Gamma(G))}{|v(\Gamma(G))|}$, where $\Gamma(G) = \Gamma_c(G)$ or $\Gamma_{nc}(G)$.
\end{theorem}
\begin{proof}
    From Theorem \ref{graphs when double toroidal}, we have that $\Gamma_c(G)$ is isomorphic to  $K_8 \sqcup 3K_4$, $K_8 \sqcup 9K_1$,  $K_8 \sqcup 5K_2$, $K_8 \sqcup 9K_3$ or $K_8 \sqcup 9(K_1 \vee 3K_2)$. If $\Gamma_c(G) \cong K_8 \sqcup 3K_4$, then $|v(\Gamma_c(G))|=20, |e(\Gamma_c(G))|=46, |e(\Gamma_{nc}(G))|= \binom{20}{2}-46=144$. Using Theorem 2.1 of \cite{SD-AS-RKN-2023} and Theorem \ref{using commuting finding complement}, we have $M_1(\Gamma_c(G))=500$, $M_2(\Gamma_c(G))=1534, M_1(\Gamma_{nc}(G))=4224$ and $M_2(\Gamma_{nc}(G))=30720$. Therefore, 
\[
\frac{M_2(\Gamma_c(G)))}{|e(\Gamma_c(G) )|}=\frac{767}{23} > 25 = \frac{M_{1}(\Gamma_c(G))}{|v(\Gamma_c(G))|}
\]
and 
\[
\frac{M_2(\Gamma_{nc}(G)))}{|e(\Gamma_{nc}(G) )|}=\frac{640}{3} > \frac{1056}{5} = \frac{M_{1}(\Gamma_{nc}(G))}{|v(\Gamma_{nc}(G))|}.
\]
If $\Gamma_c(G) \cong K_8 \sqcup 9K_1$, then $|v(\Gamma_c(G))|=17, |e(\Gamma_c(G))|=28, |e(\Gamma_{nc}(G))|= \binom{17}{2}-28=108$. Using Theorem 2.1 of \cite{SD-AS-RKN-2023} and Theorem \ref{using commuting finding complement}, we have $M_1(\Gamma_c(G))=392$, $M_2(\Gamma_c(G))=1372$, $M_1(\Gamma_{nc}(G))=2952$ and $M_2(\Gamma_{nc}(G))=19584$. Therefore, 
\[
\frac{M_2(\Gamma_c(G)))}{|e(\Gamma_c(G) )|}=49 > \frac{392}{17} = \frac{M_{1}(\Gamma_c(G))}{|v(\Gamma_c(G))|}
\]
and 
\[
\frac{M_2(\Gamma_{nc}(G)))}{|e(\Gamma_{nc}(G) )|}=\frac{1632}{9} > \frac{2952}{17} = \frac{M_{1}(\Gamma_{nc}(G))}{|v(\Gamma_{nc}(G))|}.
\]
If $\Gamma_c(G) \cong K_8 \sqcup 5K_2$, then $|v(\Gamma_c(G))|=18, |e(\Gamma_c(G))|=33, |e(\Gamma_{nc}(G))|= \binom{18}{2}-33=120$. Using Theorem 2.1 of \cite{SD-AS-RKN-2023} and Theorem \ref{using commuting finding complement}, we have $M_1(\Gamma_c(G))=402, M_2(\Gamma_c(G))=1377$, $M_1(\Gamma_{nc}(G))=3360$ and $M_2(\Gamma_{nc}(G))=23040$. Therefore, 
\[
\frac{M_2(\Gamma_c(G)))}{|e(\Gamma_c(G) )|}=\frac{67}{3} > \frac{459}{11} = \frac{M_{1}(\Gamma_c(G))}{|v(\Gamma_c(G))|}
\]
and 
\[
\frac{M_2(\Gamma_{nc}(G)))}{|e(\Gamma_{nc}(G) )|}=192 > \frac{560}{3} = \frac{M_{1}(\Gamma_{nc}(G))}{|v(\Gamma_{nc}(G))|}.
\]
If $\Gamma_c(G) \cong K_8 \sqcup 9K_3$, then $|v(\Gamma_c(G))|=35, |e(\Gamma_c(G))|=55, |e(\Gamma_{nc}(G))|= \binom{35}{2}-55=540$. Using Theorem 2.1 of \cite{SD-AS-RKN-2023} and Theorem \ref{using commuting finding complement}, we have $M_1(\Gamma_c(G))=500, M_2(\Gamma_c(G))=1480$, $M_1(\Gamma_{nc}(G))=33480$ and $M_2(\Gamma_{nc}(G))=518400$. Therefore, 
\[
\frac{M_2(\Gamma_c(G)))}{|e(\Gamma_c(G) )|}=\frac{1480}{55} > \frac{500}{35} = \frac{M_{1}(\Gamma_c(G))}{|v(\Gamma_c(G))|}
\]
and 
\[
\frac{M_2(\Gamma_{nc}(G)))}{|e(\Gamma_{nc}(G) )|}=960 > \frac{33480}{35} = \frac{M_{1}(\Gamma_{nc}(G))}{|v(\Gamma_{nc}(G))|}.
\]
If $\Gamma_c(G) \cong K_8 \sqcup 9(K_1 \vee 3K_2)$, then $|v(\Gamma_c(G))|=71, |e(\Gamma_c(G))|=109, |e(\Gamma_{nc}(G))|= \binom{71}{2}-109=2376$, $M_1(\Gamma_c(G))=932$ and $M_2(\Gamma_c(G))=2128$. Using  Theorem \ref{using commuting finding complement} we have  $M_1(\Gamma_{nc}(G))=318312$ and $M_2(\Gamma_{nc}(G))=10660608$. Therefore, 
\[
\frac{M_2(\Gamma_c(G)))}{|e(\Gamma_c(G) )|}=\frac{2128}{109} > \frac{932}{71} = \frac{M_{1}(\Gamma_c(G))}{|v(\Gamma_c(G))|}
\]
and 
\[
\frac{M_2(\Gamma_{nc}(G)))}{|e(\Gamma_{nc}(G) )|}=\frac{10660608}{2376} > \frac{318312}{71} = \frac{M_{1}(\Gamma_{nc}(G))}{|v(\Gamma_{nc}(G))|}.
\]
\end{proof}

The following theorem shows that if  $\Gamma_c(G)$ is triple-toroidal then    $\Gamma_c(G)$ and $\Gamma_{nc}(G)$ satisfy Hansen-Vuki{\v{c}}evi{\'c} Conjecture.
\begin{theorem}
    Let $G$ be a finite non-abelian group. If \, $\Gamma_c(G)$ is triple-toroidal, then $\frac{M_{2}(\Gamma(G))}{|e(\Gamma(G))|} \geq \frac{M_{1}(\Gamma(G))}{|v(\Gamma(G))|}$, where $\Gamma(G) = \Gamma_c(G)$ or $\Gamma_{nc}(G)$.
\end{theorem}
\begin{proof}
    From Theorem \ref{graphs when triple toroidal}, we have that $\Gamma_c(G)$ is isomorphic to  $6K_2 \sqcup 3K_6 \sqcup 4K_4$ or  $3K_6$. If $\Gamma_c(G) \cong 6K_2 \sqcup 3K_6 \sqcup 4K_4$, then $|v(\Gamma_c(G))|=46, |e(\Gamma_c(G))|=75, |e(\Gamma_{nc}(G))|= \binom{46}{2}-75=960$. Using Theorem 2.1 of \cite{SD-AS-RKN-2023} and Theorem \ref{using commuting finding complement}, we have $M_1(\Gamma_c(G))=606, M_2(\Gamma_c(G))=1347, M_1(\Gamma_{nc}(G))=80256$ and $M_2(\Gamma_{nc}(G))=1677120$. Therefore, 
\[
\frac{M_2(\Gamma_c(G)))}{|e(\Gamma_c(G) )|}=\frac{449}{25} > \frac{303}{23} = \frac{M_{1}(\Gamma_c(G))}{|v(\Gamma_c(G))|}
\]
and 
\[
\frac{M_2(\Gamma_{nc}(G)))}{|e(\Gamma_{nc}(G) )|}=1747 > \frac{40128}{23} = \frac{M_{1}(\Gamma_{nc}(G))}{|v(\Gamma_{nc}(G))|}.
\]
If $\Gamma_c(G) \cong 3K_6$, then $|v(\Gamma_c(G))|=18, |e(\Gamma_c(G))|=45, |e(\Gamma_{nc}(G))|= \binom{18}{2}-45=108$. Using Theorem 2.1 of \cite{SD-AS-RKN-2023} and Theorem \ref{using commuting finding complement}, we have $M_1(\Gamma_c(G))=450, M_2(\Gamma_c(G))=1125, M_1(\Gamma_{nc}(G))=2592$ and $M_2(\Gamma_{nc}(G))=15552$. Therefore, 
\[
\frac{M_2(\Gamma_c(G)))}{|e(\Gamma_c(G) )|}=25= \frac{M_{1}(\Gamma_c(G))}{|v(\Gamma_c(G))|}
\]
and 
\[
\frac{M_2(\Gamma_{nc}(G)))}{|e(\Gamma_{nc}(G) )|}= 144 = \frac{M_{1}(\Gamma_{nc}(G))}{|v(\Gamma_{nc}(G))|}.
\]
\end{proof}
Using results from \cite{PD-RKN-2021, RKN-WNF-KCD-YS-2021, MS-RKN-2023}, we have the following characterizations of $\Gamma_c(G)$ and $\Gamma_{nc}(G)$ if $\Gamma_c(G)$ is planar or toroidal.
\begin{theorem}\label{P1}
    Let $G$ be a finite non-abelian group such that  $\Gamma_c(G)$ is planar. Then
    \begin{enumerate}
        \item $\Gamma_c(G)$ is neither hyperenergetic, L-hyperenergetic nor CN-hyperenergetic.
        \item $\Gamma_c(G)$ is hypoenergetic only when $G \cong D_6$ or $D_{10}$.
        \item $\Gamma_c(G)$ is Q-hyperenergetic only when $G \cong A_4$.
        \item $E(\Gamma_c(G))<LE(\Gamma_c(G))<LE^{+}(\Gamma_c(G))$ when $G \cong A_4$ or $S_4$; $LE^+(\Gamma_c(G))<E(\Gamma_c(G))<LE(\Gamma_c(G))$ when $G \cong A_5, SL(2,3)$ or $S_z(2)$ and $E(\Gamma_c(G)) \leq LE^{+}(\Gamma_c(G)) \leq LE(\Gamma_c(G))$ otherwise. 
        \item $\Gamma_{nc}(G)$ is neither hypoenergetic nor CN-hyperenergetic.
        \item $\Gamma_{nc}(G)$ is hyperenergetic only when $G \cong S_4$.
        \item $\Gamma_{nc}(G)$ is L-hyperenergetic when $G \cong D_{6}, D_{10}, D_{12}, Q_{12}, A_{4}, A_{5}, S_{4}, SL(2,3)$ or $Sz(2)$.
        \item $\Gamma_{nc}(G)$ is Q-hyperenergetic when $G \cong D_{10}, D_{12}, Q_{12}, A_{4}, A_{5}, S_{4}$ or $SL(2,3)$.
        \item $E(\Gamma_{nc}(G)) \leq LE^{+}(\Gamma_{nc}(G)) \leq LE(\Gamma_{nc}(G))$ but $E(\Gamma_{nc}(S_4))<LE(\Gamma_{nc}(S_4))<LE^{+}(\Gamma_{nc}(S_4))$. 
    \end{enumerate}
\end{theorem}
\begin{theorem}\label{T1}
	Let $G$ be a finite non-abelian group such that  $\Gamma_c(G)$ is toroidal. Then 
	\begin{enumerate}
		\item $\Gamma_c(G)$ is neither hypoenergetc, hyperenergetic nor CN-hyperenergetic.
		\item $\Gamma_c(G)$ is L-hyperenergetic               and Q-hyperenergetic when $G \cong              D_{14}, D_{16}, Q_{16}, QD_{16}$ or             $D_6 \times \mathbb{Z}_3$.
	\item $E(\Gamma_c(G))<LE(\Gamma_c(G))<LE^+                 (\Gamma_c(G))$ when $G \cong D_6 \times             \mathbb{Z}_3$ or $A_4 \times                         \mathbb{Z}_2$ and $E(\Gamma_c(G))                  <LE^{+}(\Gamma_c(G))<LE(\Gamma_c(G))$               otherwise.
        \item $\Gamma_{nc}(G)$ is neither hypoenergetic, hyperenergetic nor CN-hyperenergetic but is L-hyper-energetic as well as Q-hyperenergetic.
        \item $E(\Gamma_{nc}(G)) \leq LE^{+}(\Gamma_{nc}(G)) \leq LE(\Gamma_{nc}(G))$ but $E(\Gamma_{nc}(A_4 \times \mathbb{Z}_2))<LE(\Gamma_{nc}(A_4 \times \mathbb{Z}_2))<LE^{+}(\Gamma_{nc}(A_4 \times \mathbb{Z}_2))$. 
    \end{enumerate}
\end{theorem}

From Theorems \ref{P1}-\ref{T1}, it follows that $\Gamma_c(G)$ and $\Gamma_{nc}(G)$ satisfy E-LE conjecture if $\Gamma_c(G)$ is planar or toroidal. In the following theorems we show that  $\Gamma_c(G)$ and $\Gamma_{nc}(G)$ satisfy E-LE conjecture if $\Gamma_c(G)$ is double-toroidal or triple-toroidal.
\begin{theorem}
	Let $G$ be a finite non-abelian group such that  $\Gamma_c(G)$ is double-toroidal. Then 
	\begin{enumerate}
		\item $\Gamma_c(G)$ is neither hyperenergetic nor CN-hyperenergetic.
		\item  $\Gamma_c(G)$ is not L-hyperenergetic only when $G \cong (\mathbb{Z}_3 \times \mathbb{Z}_3) \rtimes \mathbb{Z}_4$ or $(\mathbb{Z}_3 \times \mathbb{Z}_3) \rtimes Q_8$.
            \item $\Gamma_c(G)$ is Q-hyperenergetic.
		\item $\Gamma_c(G)$ is hypoenergetic only when $G \cong D_{18}$ or $(\mathbb{Z}_3 \times \mathbb{Z}_3) \rtimes \mathbb{Z}_{2}$.
	\item  $E(\Gamma_c(G))<LE(\Gamma_c(G))<LE^{+}(\Gamma_c(G))$ only when $G \cong (\mathbb{Z}_3 \times \mathbb{Z}_3) \rtimes \mathbb{Z}_4$ or $(\mathbb{Z}_3 \times \mathbb{Z}_3) \rtimes Q_8$ and $E(\Gamma_c(G))<LE^{+}(\Gamma_c(G))<LE(\Gamma_c(G))$ otherwise.
	\end{enumerate}  
\end{theorem}
\begin{proof}	
From Theorem \ref{graphs when double toroidal}, we have that $\Gamma_c(G)$ is isomorphic to  $K_8 \sqcup 3K_4$, $K_8 \sqcup 9K_1$,  $K_8 \sqcup 5K_2$, $K_8 \sqcup 9K_3$ or $K_8 \sqcup 9(K_1 \vee 3K_2)$. 
	
	If $\Gamma_c(G) \cong K_8 \sqcup 3K_4$, then 
	$\Spec(\Gamma_c(G))=\left \lbrace (-1)^{16}, (7)^1, (3)^3 \right \rbrace$ and so $E(\Gamma_c(G))=16+7+9=32$. We also 
	 have $\L-Spec(\Gamma_c(G))=\left \lbrace (0)^4, (8)^7, (4)^9 \right \rbrace$ and $\Q-Spec(\Gamma_c(G))=\left \lbrace (14)^1, (6)^{10}, (2)^9 \right \rbrace$. Here, $\frac{2|e(\Gamma_c(G))|}{|v(\Gamma_c(G))|}=\frac{23}{5}$ so $|0-\frac{23}{5}|=\frac{23}{5}, |8-\frac{23}{5}|=\frac{17}{5}, |4-\frac{23}{5}|=\frac{3}{5}$. Therefore, $LE(\Gamma_c(G))=4 \cdot \frac{23}{5}+7 \cdot \frac{17}{5}+9 \cdot \frac{3}{5}=\frac{238}{5}$. Similarly, $|14-\frac{23}{5}|=\frac{47}{5}, |6-\frac{23}{5}|=\frac{7}{5}, |2-\frac{23}{5}|=\frac{13}{5}$ and hence $LE^{+}(\Gamma_c(G))=\frac{47}{5}+10 \cdot \frac{7}{5}+9 \cdot \frac{13}{5}=\frac{234}{5}$. 
	 Further, $\CN-Spec(\Gamma_c(G))=\left \lbrace (-6)^7, (42)^1, (-2)^9, (6)^3 \right \rbrace$ and so $E_{CN}(\Gamma_c(G))=120$. Now, $|v(\Gamma_c(G))|=20< 32=E(\Gamma_c(G)), E(K_{20})=2(20-1)=38 > 32=E(\Gamma_c(G))$ and $E_{CN}(K_{20})=2(20-1)(20-2)=684>120=E_{CN}(\Gamma_c(G))$. Thus, \textit{$\Gamma_c(G)$ is neither hypoenergetic, hyperenergetic  nor CN-hyperenergetic}. Also, $LE(K_{20})=2(20-1)=38<\frac{238}{5}=LE(\Gamma_c(G))$ and $LE^{+}(K_{20})=2(20-1)=38<\frac{234}{5}=LE^{+}(\Gamma_c(G))$. Therefore, \textit{$\Gamma_c(G)$ is L-hyperenergetic as well as Q-hyperenergetic}. Further, 
	 \[
	 E(\Gamma_c(G))=32<\frac{234}{5}=LE^{+}(\Gamma_c(G))<\frac{238}{5}=LE(\Gamma_c(G)).
	 \]
	
	If $\Gamma_c(G) \cong K_8 \sqcup 9K_1$, then 
	$\Spec(\Gamma_c(G))=\left \lbrace (-1)^{7}, (7)^1, (0)^9 \right \rbrace$ and so $E(\Gamma_c(G))=7+7=14$. 
	We also have $\L-Spec(\Gamma_c(G))=\left \lbrace (0)^{10}, (8)^7, \right \rbrace$ and $\Q-Spec(\Gamma_c(G))=\left \lbrace (14)^1, (6)^{7}, (0)^9 \right \rbrace$. Here, $\frac{2|e(\Gamma_c(G))|}{|v(\Gamma_c(G))|}=\frac{56}{17}$ so $|0-\frac{56}{17}|=\frac{56}{17}, |8-\frac{56}{17}|=\frac{80}{17}$. Therefore, $LE(\Gamma_c(G))=10 \cdot \frac{56}{17}+7 \cdot \frac{80}{17}=\frac{1120}{17}$. Similarly, $|14-\frac{56}{17}|=\frac{182}{17}, |6-\frac{56}{17}|=\frac{46}{17}, |0-\frac{56}{17}|=\frac{56}{17}$ and hence $LE^{+}(\Gamma_c(G))=\frac{182}{17}+7 \cdot \frac{46}{17}+9 \cdot \frac{56}{17}=\frac{1008}{17}$. 
	Further, $\CN-Spec(\Gamma_c(G))=\left \lbrace (-6)^7, (42)^1, (0)^9 \right \rbrace$ and so $E_{CN}(\Gamma_c(G))=84$. Now, $|v(\Gamma_c(G))|=17>14=E(\Gamma_c(G)), E(K_{17})=2(17-1)=32 > 14=E(\Gamma_c(G))$ and $E_{CN}(K_{17})=2(17-1)(17-2)=480>84=E_{CN}(\Gamma_c(G))$. Thus, \textit{$\Gamma_c(G)$ is hypoenergetic but neither hyperenergetic nor CN-hyperenergetic}. Also, $LE(K_{17})=2(17-1)=32<\frac{1120}{17}=LE(\Gamma_c(G))$ and $LE^{+}(K_{17})=2(17-1)=32<\frac{1008}{17}=LE^{+}(\Gamma_c(G))$. Therefore, \textit{$\Gamma_c(G)$ is L-hyperenergetic as well as Q-hyperenergetic}. Further, 
	\[
	E(\Gamma_c(G))=14<\frac{1008}{17}=LE^{+}(\Gamma_c(G))<\frac{1120}{17}=LE(\Gamma_c(G)).
	\]
	
	If $\Gamma_c(G) \cong K_8 \sqcup 5K_2$, then 
	$\Spec(\Gamma_c(G))=\left \lbrace (-1)^{12}, (7)^1, (1)^5 \right \rbrace$ and so $E(\Gamma_c(G))=12+7+5=24$. 
	We also have $\L-Spec(\Gamma_c(G))=\left \lbrace (0)^6, (8)^7, (2)^5 \right \rbrace$ and $\Q-Spec(\Gamma_c(G))=\left \lbrace (14)^1, (6)^{7}, (2)^5, (0)^5 \right \rbrace$. Here, $\frac{2|e(\Gamma_c(G))|}{|v(\Gamma_c(G))|}=\frac{11}{3}$ so $|0-\frac{11}{3}|=\frac{11}{3}, |8-\frac{11}{3}|=\frac{13}{3}, |2-\frac{11}{3}|=\frac{5}{3}$. Therefore, $LE(\Gamma_c(G))=6 \cdot \frac{11}{3}+7 \cdot \frac{13}{3}+5 \cdot \frac{5}{3}=\frac{182}{3}$. Similarly, $|14-\frac{11}{3}|=\frac{31}{3}, |6-\frac{11}{3}|=\frac{7}{3}, |2-\frac{11}{3}|=\frac{5}{3}, |0-\frac{11}{3}|=\frac{11}{3}$ and hence $LE^{+}(\Gamma_c(G))=\frac{31}{3}+7 \cdot \frac{7}{3}+5 \cdot \frac{5}{3}+ 5 \cdot \frac{11}{3}=\frac{160}{3}$. 
	Further, $\CN-Spec(\Gamma_c(G))=\left \lbrace (-6)^7, (42)^1, (0)^{10} \right \rbrace$ and so $E_{CN}(\Gamma_c(G))=84$. Now, $|v(\Gamma_c(G))|=18< 24=E(\Gamma_c(G)), E(K_{18})=2(18-1)=34 > 24=E(\Gamma_c(G))$ and $E_{CN}(K_{18})=2(18-1)(18-2)=544>84=E_{CN}(\Gamma_c(G))$. Thus, \textit{$\Gamma_c(G)$ is neither hypoenergetic,  hyperenergetic  nor CN-hyperenergetic}. Also, $LE(K_{18})=2(18-1)=34<\frac{182}{3}=LE(\Gamma_c(G))$ and $LE^{+}(K_{18})=2(18-1)=34<\frac{160}{3}=LE^{+}(\Gamma_c(G))$. Therefore, \textit{$\Gamma_c(G)$ is L-hyperenergetic as well as Q-hyperenergetic}. Further,
	\[
	E(\Gamma_c(G))=24<\frac{160}{3}=LE^{+}(\Gamma_c(G))<\frac{182}{3}=LE(\Gamma_c(G)).
	\]
	If $\Gamma_c(G) \cong K_8 \sqcup 9K_3$, then 
	$\Spec(\Gamma_c(G))=\left \lbrace (-1)^{25}, (7)^1, (2)^9 \right \rbrace$ and so $E(\Gamma_c(G))=25+7+18=50$. 
	We also have $\L-Spec(\Gamma_c(G))=\left \lbrace (0)^{10}, (8)^7, (3)^{18} \right \rbrace$ and $\Q-Spec(\Gamma_c(G))=\left \lbrace (14)^1, (6)^{7}, (4)^9, (1)^{18} \right \rbrace$. Here, $\frac{2|e(\Gamma_c(G))|}{|v(\Gamma_c(G))|}=\frac{22}{7}$ so $|0-\frac{22}{7}|=\frac{22}{7}, |8-\frac{22}{7}|=\frac{34}{7}, |3-\frac{22}{7}|=\frac{1}{7}$. Therefore, $LE(\Gamma_c(G))=10 \cdot \frac{22}{7}+7 \cdot \frac{34}{7}+18 \cdot \frac{1}{7}=68$. Similarly, $|14-\frac{22}{7}|=\frac{76}{7}, |6-\frac{22}{7}|=\frac{20}{7}, |4-\frac{22}{7}|=\frac{6}{7}, |1-\frac{22}{7}|=\frac{15}{7}$ and hence $LE^{+}(\Gamma_c(G))=\frac{76}{7}+7 \cdot \frac{20}{7}+9 \cdot \frac{6}{7}+ 18 \cdot \frac{15}{7}=\frac{540}{7}$. 
	Further, $\CN-Spec(\Gamma_c(G))=\left \lbrace (-6)^7, (42)^1, (-1)^{18}, (2)^9 \right \rbrace$ and so $E_{CN}(\Gamma_c(G))=120$. Now, $|v(\Gamma_c(G))|=35< 50=E(\Gamma_c(G)),$ $ E(K_{35})=2(35-1)=68 > 50=E(\Gamma_c(G))$ and $E_{CN}(K_{35})=2(35-1)(35-2)=2244>120=E_{CN}(\Gamma_c(G))$. Thus, \textit{$\Gamma_c(G)$ is neither hypoenergetic, hyperenergetic  nor CN-hyperenergetic}. Also, $LE(K_{35})=2(35-1)=68=LE(\Gamma_c(G))$ and $LE^{+}(K_{35})=2(35-1)=68<\frac{540}{7}=LE^{+}(\Gamma_c(G))$. Therefore, \textit{$\Gamma_c(G)$ is not L-hyperenergetic but Q-hyperenergetic}. Further,
	\[
	E(\Gamma_c(G))=50<68=LE(\Gamma_c(G))<\frac{540}{7}=LE^{+}(\Gamma_c(G)).
	\]
 If $\Gamma_c(G) \cong K_8 \sqcup 9(K_1 \vee 3K_2)$, then 
	$\Spec(\Gamma_c(G))=\left \lbrace (-1)^{34}, (7)^1, (-2)^9, (1)^{18}, (3)^9 \right \rbrace$ and so $E(\Gamma_c(G))=34+34+18+18=104$. 
	We also have $\L-Spec(\Gamma_c(G))=\left \lbrace (0)^{10}, (8)^7, (3)^{27}, (1)^{18}, (7)^{9} \right \rbrace$ and $\Q-Spec(\Gamma_c(G))=\left \lbrace (14)^1, (6)^{7}, (3)^{18}, (1)^{27}, \left(\frac{9+\sqrt{33}}{2} \right)^9, \left(\frac{9-\sqrt{33}}{2} \right)^9 \right \rbrace$. Here, $\frac{2|e(\Gamma_c(G))|}{|v(\Gamma_c(G))|}=\frac{218}{71}$ so $|0-\frac{218}{71}|=\frac{218}{71}, |8-\frac{218}{71}|=\frac{350}{71}, |3-\frac{218}{71}|=\frac{5}{71}, |1-\frac{218}{71}|=\frac{147}{71}, |7-\frac{218}{71}|=\frac{279}{71}$. Therefore, $LE(\Gamma_c(G))=10 \cdot \frac{218}{71}+7 \cdot \frac{350}{71}+18 \cdot \frac{147}{71}+27 \cdot \frac{5}{71}+9 \cdot \frac{279}{71}=\frac{9922}{71}$. Similarly, $|14-\frac{218}{71}|=\frac{776}{71}, |6-\frac{218}{71}|=\frac{208}{71}, |\frac{9+\sqrt{33}}{2}-\frac{218}{71}| \approx \frac{610.86}{71}, |\frac{9-\sqrt{33}}{2}-\frac{218}{71}| \approx \frac{204.86}{71}$ and hence $LE^{+}(\Gamma_c(G)) \approx \frac{776}{71}+7 \cdot \frac{208}{71}+27 \cdot \frac{147}{71}+18 \cdot \frac{5}{71}+9 \cdot \frac{610.86}{71}+ 9 \cdot \frac{204.86}{71} \approx \frac{13632.48}{71}$. 
	Further, $\CN-Spec(\Gamma_c(G))=\left \lbrace (-6)^7, (42)^1, (-1)^{54}, (6)^9 \right \rbrace$ and so $E_{CN}(\Gamma_c(G))=192$. Now, $|v(\Gamma_c(G))|=71< 104=E(\Gamma_c(G)), E(K_{71})=2(71-1)=140 > 104=E(\Gamma_c(G))$ and $E_{CN}(K_{71})=2(71-1)(71-2)=9660>192=E_{CN}(\Gamma_c(G))$. Thus, \textit{$\Gamma_c(G)$ is neither hypoenergetic, hyperenergetic  nor CN-hyperenergetic}. Also, $LE(K_{71})=2(71-1)=140 > \frac{9922}{71}=LE(\Gamma_c(G))$ and $LE^{+}(K_{71})=2(71-1)=140<\frac{13632.48}{71} \approx LE^{+}(\Gamma_c(G))$. Therefore, \textit{$\Gamma_c(G)$ is not L-hyperenergetic but Q-hyperenergetic}. Further,
	\[
	E(\Gamma_c(G))=104<\frac{9922}{71}=LE(\Gamma_c(G))<\frac{13632.48}{71}=LE^{+}(\Gamma_c(G)).
	\]
\end{proof}
\begin{theorem}
	Let $G$ be a finite non-abelian group such that $\Gamma_c(G)$ is triple-toroidal. Then 
	\begin{enumerate}
		\item $\Gamma_c(G)$ is neither hypoenergetic, hyperenergetic, CN-hyperenergetic  nor Q-hyperenergetic.
		\item  $\Gamma_c(G)$ is L-hyperenergetic only when $G \cong GL(2, 3)$. 
		\item $E(\Gamma_c(G)) \leq LE^{+}(\Gamma_c(G)) \leq LE(\Gamma_c(G))$.
	\end{enumerate}
\end{theorem}
\begin{proof}
	From Theorem \ref{graphs when triple toroidal}, we have that $\Gamma_c(G)$ is isomorphic to  $6K_2 \sqcup 3K_6 \sqcup 4K_4$ or $3K_6$.
	
	If $\Gamma_c(G) \cong 6K_2 \sqcup 3K_6 \sqcup 4K_4$, then 
	$\Spec(\Gamma_c(G))=\left \lbrace (-1)^{33}, (1)^6, (5)^3, (3)^4 \right \rbrace$ and so $E(\Gamma_c(G))=33+6+15+12=66$. 
	We also have $\L-Spec(\Gamma_c(G))=\left \lbrace (0)^{13}, (2)^6, (6)^{15}, (4)^{12} \right \rbrace$ and $\Q-Spec(\Gamma_c(G))=\left \lbrace (0)^6, (10)^3, (4)^{15}, (6)^4, (2)^{18} \right \rbrace$. Here, $\frac{2|e(\Gamma_c(G))|}{|v(\Gamma_c(G))|}=\frac{75}{23}$ so $|0-\frac{75}{23}|=\frac{75}{23}, |2-\frac{75}{23}|=\frac{29}{23}, |6-\frac{75}{23}|=\frac{63}{23}, |4-\frac{75}{23}|=\frac{17}{23}$. Therefore, $LE(\Gamma_c(G))=13 \cdot \frac{75}{23}+6 \cdot \frac{29}{23}+15 \cdot \frac{63}{23}+ 12 \cdot \frac{17}{23}=\frac{2298}{23}$. Similarly, $|0-\frac{75}{23}|=\frac{75}{23}, |10-\frac{75}{23}|=\frac{155}{23}, |4-\frac{75}{23}|=\frac{17}{23}, |6-\frac{75}{23}|=\frac{63}{23}, |2-\frac{75}{23}|=\frac{29}{23}$ and hence $LE^{+}(\Gamma_c(G))=6 \cdot \frac{75}{23}+3 \cdot \frac{155}{23}+15 \cdot\frac{17}{23}+4 \cdot \frac{63}{23}+ 18 \cdot \frac{29}{23}=\frac{1944}{23}$. 
	Further, $\CN-Spec(\Gamma_c(G))=\left \lbrace (0)^{12}, (-4)^{15}, (20)^3, (-2)^{12}, (6)^4 \right \rbrace$ and so $E_{CN}(\Gamma_c(G))=168$. Now, $|v(\Gamma_c(G))|=46< 66=E(\Gamma_c(G)), E(K_{46})=2(46-1)=90 > 66=E(\Gamma_c(G))$ and $E_{CN}(K_{46})=2(46-1)(46-2)=3960>168=E_{CN}(\Gamma_c(G))$. Thus, \textit{$\Gamma_c(G)$ is neither hypoenergetic, hyperenergetic nor CN-hyperenergetic}. Also, $LE(K_{46})=2(46-1)=90<\frac{2298}{23}=LE(\Gamma_c(G))$ and $LE^{+}(K_{46})=2(46-1)=90>\frac{1944}{23}=LE^{+}(\Gamma_c(G))$. Therefore, \textit{$\Gamma_c(G)$ is L-hyperenergetic but not Q-hyperenergetic}. Further, 
	\[
	E(\Gamma_c(G))=66<\frac{1944}{23}=LE^{+}(\Gamma_c(G))<\frac{2298}{23}=LE(\Gamma_c(G)).
	\]
	
	If $\Gamma_c(G) \cong 3K_6$, then 
	$\Spec(\Gamma_c(G))=\left \lbrace (-1)^{15}, (5)^3 \right \rbrace$ and so $E(\Gamma_c(G))=15+15=30$. 
	We also have $\L-Spec(\Gamma_c(G))=\left \lbrace (0)^{3}, (6)^{15} \right \rbrace$ and $\Q-Spec(\Gamma_c(G))=\left \lbrace (10)^3, (4)^{15} \right \rbrace$. Here, $\frac{2|e(\Gamma_c(G))|}{|v(\Gamma_c(G))|}=5$ so $|0-5|=5, |6-5|=1$. Therefore, $LE(\Gamma_c(G))=3\cdot 5 + 15 \cdot 1=30$. Similarly, $|10-5|=5, |4-5|=1$ and hence $LE^{+}(\Gamma_c(G))=3 \cdot 5 + 15 \cdot 1=30$. 
	Further, $\CN-Spec(\Gamma_c(G))=\left \lbrace (-4)^{15}, (20)^3 \right \rbrace$ and so $E_{CN}(\Gamma_c(G))=120$. Now, $|v(\Gamma_c(G))|=18< 30=E(\Gamma_c(G)), E(K_{18})=2(18-1)=34 > 30=E(\Gamma_c(G))$ and $E_{CN}(K_{18})=2(18-1)(18-2)=544>120=E_{CN}(\Gamma_c(G))$. Thus, \textit{$\Gamma_c(G)$ is neither hypoenergetic, hyperenergetic  nor CN-hyperenergetic}. Also, $LE(K_{18})=2(18-1)=34>30=LE(\Gamma_c(G))$ and $LE^{+}(K_{18})=2(18-1)=34>30=LE^{+}(\Gamma_c(G))$. Therefore, \textit{$\Gamma_c(G)$ is neither L-hyperenergetic nor Q-hyperenergetic}. Further,
	\[
	E(\Gamma_c(G))=30=LE^{+}(\Gamma_c(G))=LE(\Gamma_c(G)).
	\]
\end{proof}
\begin{theorem}
	Let $G$ be a finite non-abelian group such that  $\Gamma_c(G)$ is double-toroidal. Then 
	\begin{enumerate}
		\item $\Gamma_{nc}(G)$ is neither hypoenergetic nor CN-hyperenergetic.
            \item $\Gamma_{nc}(G)$ is hyperenergetic only when $G \cong (\mathbb{Z}_3 \times \mathbb{Z}_3) \rtimes Q_8$.
		\item  $\Gamma_{nc}(G)$ is L-hyperenergetic and Q-hyperenergetic.
	\item  $E(\Gamma_{nc}(G))<LE(\Gamma_{nc}(G))<LE^{+}(\Gamma_{nc}(G))$ only when $G \cong (\mathbb{Z}_3 \times \mathbb{Z}_3) \rtimes \mathbb{Z}_4$ and $E(\Gamma_{nc}(G))<LE^{+}(\Gamma_{nc}(G))<LE(\Gamma_{nc}(G))$ otherwise.
	\end{enumerate}  
\end{theorem}
\begin{proof}	
From Theorem \ref{graphs when double toroidal}, we have that $\Gamma_c(G)$ is isomorphic to  $K_8 \sqcup 3K_4$, $K_8 \sqcup 9K_1$,  $K_8 \sqcup 5K_2$, $K_8 \sqcup 9K_3$ or $K_8 \sqcup 9(K_1 \vee 3K_2)$.
	
	If $\Gamma_c(G) \cong K_8 \sqcup 3K_4$, then 
	$\Spec(\Gamma_{nc}(G))=\left \lbrace (0)^{16}, (-4)^2, (4+\sqrt{112})^1, (4-\sqrt{112})^1 \right \rbrace$ and so $E(\Gamma_{nc}(G))=8+2\sqrt{112}$. We also 
	 have $\L-Spec(\Gamma_{nc}(G))=\left \lbrace (0)^1, (16)^9, (12)^7, (20)^3 \right \rbrace$ and $\Q-Spec(\Gamma_{nc}(G))=\left \lbrace (12)^9, (16)^{9}, (18+\sqrt{132})^1, (18-\sqrt{132})^1 \right \rbrace$. Here, $\frac{2|e(\Gamma_{nc}(G))|}{|v(\Gamma_{nc}(G))|}=\frac{72}{5}$ so $|0-\frac{72}{5}|=\frac{72}{5}, |16-\frac{72}{5}|=\frac{8}{5}, |12-\frac{72}{5}|=\frac{12}{5}, |20-\frac{72}{5}|=\frac{28}{5}$. Therefore, $LE(\Gamma_{nc}(G))= \frac{72}{5}+9 \cdot \frac{8}{5}+7 \cdot \frac{12}{5}+3 \cdot \frac{28}{5}=\frac{312}{5}$. Similarly, $|12-\frac{72}{5}|=\frac{12}{5}, |16-\frac{72}{5}|=\frac{7}{5}, |18+\sqrt{132}-\frac{72}{5}|=\frac{18+5\sqrt{132}}{5}, |18-\sqrt{132}-\frac{72}{5}|=\frac{5\sqrt{132}-18}{5}$ and hence $LE^{+}(\Gamma_{nc}(G))=9 \cdot \frac{12}{5}+9 \cdot \frac{8}{5}+ \frac{18+5\sqrt{132}}{5}+\frac{5\sqrt{132}-18}{5}=36+2\sqrt{132}$. 
	 Further, $\CN-Spec(\Gamma_{nc}(G))=\left \lbrace 2(57+\sqrt{1761})^1, 2(57-\sqrt{1761})^1, (-16)^9, (-12)^7, (0)^2 \right \rbrace$ and so $E_{CN}(\Gamma_{nc}(G))=456$. Now, $|v(\Gamma_{nc}(G))|=20< 8+2\sqrt{112}=E(\Gamma_{nc}(G)), E(K_{20})=2(20-1)=38 > 8+2\sqrt{112}=E(\Gamma_{nc}(G))$ and $E_{CN}(K_{20})=2(20-1)(20-2)=684>456=E_{CN}(\Gamma_{nc}(G))$. Thus, \textit{$\Gamma_{nc}(G)$ is neither hypoenergetic, hyperenergetic  nor CN-hyperenergetic}. Also, $LE(K_{20})=2(20-1)=38<\frac{312}{5}=LE(\Gamma_{nc}(G))$ and $LE^{+}(K_{20})=2(20-1)=38<36+2\sqrt{132}=LE^{+}(\Gamma_{nc}(G))$. Therefore, \textit{$\Gamma_{nc}(G)$ is L-hyperenergetic as well as Q-hyperenergetic}. Further, 
	 \[
	 E(\Gamma_{nc}(G))=8+2\sqrt{112}<36+2\sqrt{132}=LE^{+}(\Gamma_{nc}(G))<\frac{312}{5}=LE(\Gamma_{nc}(G)).
	 \]
	
	If $\Gamma_c(G) \cong K_8 \sqcup 9K_1$, then 
	$\Spec(\Gamma_{nc}(G))=\left \lbrace (0)^{7}, (-1)^8, (4+\sqrt{88})^1, (4-\sqrt{88})^1 \right \rbrace$ and so $E(\Gamma_{nc}(G))=8+2\sqrt{88}$. We also 
	 have $\L-Spec(\Gamma_{nc}(G))=\left \lbrace (0)^1, (9)^7, (17)^9 \right \rbrace$ and $\Q-Spec(\Gamma_{nc}(G))=\left \lbrace (9)^7, (15)^{8}, \left(\frac{33+\sqrt{513}}{2} \right)^1, \left(\frac{33-\sqrt{513}}{2} \right)^1 \right \rbrace$. Here, $\frac{2|e(\Gamma_{nc}(G))|}{|v(\Gamma_{nc}(G))|}=\frac{216}{17}$ so $|0-\frac{216}{17}|=\frac{216}{17}, |9-\frac{216}{17}|=\frac{63}{17}, |17-\frac{216}{17}|=\frac{73}{17}$. Therefore, $LE(\Gamma_{nc}(G))= \frac{216}{17}+7\cdot \frac{63}{17}+9 \cdot \frac{73}{17}=\frac{1314}{17}$. Similarly, $|9-\frac{216}{17}|=\frac{63}{17}, |15-\frac{216}{17}|=\frac{39}{17}, |\frac{33+\sqrt{513}}{2}-\frac{216}{17}|=\frac{129+17\sqrt{513}}{34}, |\frac{33-\sqrt{513}}{2}-\frac{216}{17}|=\frac{17\sqrt{513}-129}{34}$ and hence $LE^{+}(\Gamma_{nc}(G))=7 \cdot \frac{63}{17}+8 \cdot \frac{39}{17}+ \frac{129+17\sqrt{513}}{34}+\frac{17\sqrt{513}-129}{34}=\frac{753+17\sqrt{513}}{17}$. 
	 Further, $\CN-Spec(\Gamma_{nc}(G))=\left \lbrace \frac{3}{2}(61+\sqrt{2049})^1, \right .$ $ \left . \frac{3}{2}(61-\sqrt{2049})^1, (-15)^8, (-9)^7 \right \rbrace$ and so $E_{CN}(\Gamma_{nc}(G))=366$. Now, $|v(\Gamma_{nc}(G))|=17< 8+2\sqrt{88}=E(\Gamma_{nc}(G)), E(K_{17})=2(17-1)=32 > 8+2\sqrt{88}=E(\Gamma_{nc}(G))$ and $E_{CN}(K_{17})=2(17-1)(17-2)=480>366=E_{CN}(\Gamma_{nc}(G))$. Thus, \textit{$\Gamma_{nc}(G)$ is neither hypoenergetic, hyperenergetic  nor CN-hyperenergetic}. Also, $LE(K_{17})=2(17-1)=32<\frac{1314}{17}=LE(\Gamma_{nc}(G))$ and $LE^{+}(K_{17})=2(17-1)=32<\frac{753+17\sqrt{513}}{17}=LE^{+}(\Gamma_{nc}(G))$. Therefore, \textit{$\Gamma_{nc}(G)$ is L-hyperenergetic as well as Q-hyperenergetic}. Further, 
	 \[
	 E(\Gamma_{nc}(G))=8+2\sqrt{88}<\frac{753+17\sqrt{513}}{17}=LE^{+}(\Gamma_{nc}(G))<\frac{1314}{17}=LE(\Gamma_{nc}(G)).
	 \]
	
	If $\Gamma_c(G) \cong K_8 \sqcup 5K_2$, then 
	$\Spec(\Gamma_{nc}(G))=\left \lbrace (0)^{12}, (-2)^6, (4+\sqrt{96})^1, (4-\sqrt{96})^1 \right \rbrace$ and so $E(\Gamma_{nc}(G))=12+2\sqrt{96}$. We also 
	 have $\L-Spec(\Gamma_{nc}(G))=\left \lbrace (0)^1, (16)^5, (10)^7, (18)^5 \right \rbrace$ and \\ $\Q-Spec(\Gamma_{nc}(G))=\left \lbrace (10)^7, (16)^{5}, (14)^4, (17+\sqrt{129})^1, (17-\sqrt{129})^1 \right \rbrace$. Here, $\frac{2|e(\Gamma_{nc}(G))|}{|v(\Gamma_{nc}(G))|}=\frac{40}{3}$ so $|0-\frac{40}{3}|=\frac{40}{3}, |16-\frac{40}{3}|=\frac{8}{3}, |10-\frac{40}{3}|=\frac{10}{3}, |18-\frac{40}{3}|=\frac{14}{3}$. Therefore, $LE(\Gamma_{nc}(G))= \frac{40}{3}+5 \cdot \frac{8}{3}+7 \cdot \frac{10}{3}+5 \cdot \frac{14}{3}=\frac{220}{3}$. Similarly, $|10-\frac{40}{3}|=\frac{10}{3}, |16-\frac{40}{3}|=\frac{8}{3}, |14-\frac{40}{3}|=\frac{2}{3},  |17+\sqrt{129}-\frac{40}{3}|=\frac{11+3\sqrt{129}}{3}, |17-\sqrt{129}-\frac{40}{3}|=\frac{3\sqrt{129}-11}{3}$ and hence $LE^{+}(\Gamma_{nc}(G))=7 \cdot \frac{10}{3}+5 \cdot \frac{8}{3}+ 4 \cdot \frac{2}{3}+ \frac{11+3\sqrt{129}}{3}+\frac{3\sqrt{129}-11}{3}=\frac{118+6\sqrt{129}}{3}$. 
	 Further, $\CN-Spec(\Gamma_{nc}(G))=\left \lbrace (99+\sqrt{5961})^1, (99-\sqrt{5961})^1, (-16)^5, (-2)^4, (-10)^7 \right \rbrace$ and so $E_{CN}(\Gamma_{nc}(G))=356$. Now, $|v(\Gamma_{nc}(G))|=18< 12+2\sqrt{96}=E(\Gamma_{nc}(G)), E(K_{18})=2(18-1)=34 > 12+2\sqrt{96}=E(\Gamma_{nc}(G))$ and $E_{CN}(K_{18})=2(18-1)(18-2)=544>356=E_{CN}(\Gamma_{nc}(G))$. Thus, \textit{$\Gamma_{nc}(G)$ is neither hypoenergetic, hyperenergetic  nor CN-hyperenergetic}. Also, $LE(K_{18})=2(18-1)=34<\frac{220}{3}=LE(\Gamma_{nc}(G))$ and $LE^{+}(K_{18})=2(18-1)=34<\frac{118+6\sqrt{129}}{3}=LE^{+}(\Gamma_{nc}(G))$. Therefore, \textit{$\Gamma_{nc}(G)$ is L-hyperenergetic as well as Q-hyperenergetic}. Further, 
	 \[
	 E(\Gamma_{nc}(G))=12+2\sqrt{96}<\frac{118+6\sqrt{129}}{3}=LE^{+}(\Gamma_{nc}(G))<\frac{220}{3}=LE(\Gamma_{nc}(G)).
	 \]
If $\Gamma_c(G) \cong K_8 \sqcup 9K_3$, then 
	$\Spec(\Gamma_{nc}(G))=\left \lbrace (0)^{25}, (-3)^8, \left(12+6\sqrt{10} \right)^1, \left(12-6\sqrt{10}\right)^1 \right \rbrace$ and so $E(\Gamma_{nc}(G))=24+12\sqrt{10}$. We also have $\L-Spec(\Gamma_{nc}(G))=\left \lbrace (0)^1, (27)^7, (32)^{18}, (35)^9 \right \rbrace$ and \\ $\Q-Spec(\Gamma_{nc}(G))=\left \lbrace (27)^7, (29)^{8}, (32)^{18}, \left(\frac{83+\sqrt{12073}}{2}\right)^1, \left(\frac{83-\sqrt{12073}}{2}\right)^1 \right \rbrace$. Here, $\frac{2|e(\Gamma_{nc}(G))|}{|v(\Gamma_{nc}(G))|}=\frac{216}{7}$ so $|0-\frac{216}{7}|=\frac{216}{7}, |32-\frac{216}{7}|=\frac{8}{7}, |27-\frac{216}{7}|=\frac{27}{7}, |35-\frac{216}{7}|=\frac{29}{7}$. Therefore, $LE(\Gamma_{nc}(G))= \frac{216}{7}+18 \cdot \frac{8}{7}+7 \cdot \frac{27}{7}+9 \cdot \frac{29}{7}=\frac{810}{7}$. Similarly, $|29-\frac{216}{7}|=\frac{13}{7}, |\frac{83+\sqrt{12073}}{2}-\frac{216}{7}| \approx \frac{918.14}{14}, |\frac{83-\sqrt{12073}}{2}-\frac{216}{7}| \approx \frac{620.14}{14}$ and hence $LE^{+}(\Gamma_{nc}(G))=7 \cdot \frac{27}{7}+18 \cdot \frac{8}{7}+ 8 \cdot \frac{13}{7}+ \frac{918.14}{14}+\frac{620.14}{14} \approx \frac{2412.28}{14}$. 
	 Further, $\CN-Spec(\Gamma_{nc}(G))=\left \lbrace \left(\frac{949+\sqrt{823705}}{2}\right)^1, \left(\frac{949-\sqrt{823705}}{2}\right)^1, (-32)^{18}, (-27)^7, (-23)^8 \right \rbrace$ and so $E_{CN}(\Gamma_{nc}(G))=1898$. Now, \\
     $|v(\Gamma_{nc}(G))|=35< 24+12\sqrt{10}=E(\Gamma_{nc}(G)), E(K_{35})=2(35-1)=68 > 24+12\sqrt{10}=E(\Gamma_{nc}(G))$ and $E_{CN}(K_{35})=2(35-1)(35-2)=2244>1898=E_{CN}(\Gamma_{nc}(G))$. Thus, \textit{$\Gamma_{nc}(G)$ is neither hypoenergetic, hyperenergetic  nor CN-hyperenergetic}. Also, $LE(K_{35})=2(35-1)=68<\frac{810}{7}=LE(\Gamma_{nc}(G))$ and $LE^{+}(K_{35})=2(35-1)=68<\frac{2412.28}{14}=LE^{+}(\Gamma_{nc}(G))$. Therefore, \textit{$\Gamma_{nc}(G)$ is L-hyperenergetic as well as Q-hyperenergetic}. Further, 
	 \[
	 E(\Gamma_{nc}(G))=24+12\sqrt{10}<\frac{810}{7}=LE(\Gamma_{nc}(G))<\frac{2412.28}{14}=LE^{+}(\Gamma_{nc}(G)).
	 \]
If $\Gamma_c(G) \cong K_8 \sqcup 9(K_1 \vee 3K_2)$, then $\Spec(\Gamma_{nc}(G))=\left \lbrace (0)^{34}, (-2)^{18}, (-4)^8, (1)^8, (x_1)^1, \right .$
 $\left . (x_2)^1, (x_3)^1 \right \rbrace$, where $x_1, x_2$ and $x_3$ are roots of the equation $x^3-60x^2-472x+288=0$.
 Since  $x_1 \approx 66.98, x_2 \approx -7.55, x_3 \approx 0.569$, we have $E(\Gamma_{nc}(G))= 8+36+32+66.98+7.55+0.569 \approx 151.09$. 
	We also have $\L-Spec(\Gamma_{nc}(G))=\left \lbrace (71)^{7}, (70)^{16}, (68)^{27}, (64)^{7}, (63)^{7}, (y_1)^1, (y_2)^1, (y_3)^1, (z_1)^1, (z_2)^1, (z_3)^1, (z_4)^1 \right \rbrace$, where $y_1, y_2$ and $y_3$ are roots of the equation $x^3-205x^2+13994x-318088$ and $z_1, z_2, z_3$ and $z_4$ are roots of the equation $x^4-205x^3+14010x^2-320232x+71680$ and $\Q-Spec(\Gamma_{nc}(G))$ $=\left \lbrace (68)^{27}, (66)^{18}, (63)^{7}, \left(\frac{129+\sqrt{33}}{2}\right)^8, \left(\frac{129-\sqrt{33}}{2}\right)^{8}, (l_1)^1, (l_2)^1, (l_3)^1 \right \rbrace$, where $l_1, l_2$ and $l_3$ are roots of the equation $x^3-255x^2+19848x-487296=0$. Here, $\frac{2|e(\Gamma_{nc}(G))|}{|v(\Gamma_{nc}(G))|}=\frac{4752}{71}$ so $|71-\frac{4752}{71}|=\frac{289}{71}, |70-\frac{4752}{71}|=\frac{218}{71}, |68-\frac{4752}{71}|=\frac{76}{71}, |64-\frac{4752}{71}|=\frac{208}{71}, |63-\frac{4752}{71}|=\frac{279}{71}$.Since $y_1 \approx 71.63, y_2 \approx 69.07, y_3 \approx 64.20, z_1 \approx 71.49, z_2 \approx 69.15, z_3 \approx 64.21$ and $z_4 \approx 0.226$, we have $|y_1-\frac{4752}{71}| \approx \frac{333.73}{71}, |y_2-\frac{4752}{71}| \approx \frac{151.97}{71}, |y_3-\frac{4752}{71}| \approx \frac{193.8}{71}, |z_1-\frac{4752}{71}| \approx \frac{323.79}{71}, |z_2-\frac{4752}{71}| \approx \frac{157.65}{71}, |z_3-\frac{4752}{71}| \approx \frac{193.09}{71}$ and $|z_4-\frac{4752}{71}| \approx \frac{4736.38}{71}$. Therefore, $LE(\Gamma_{nc}(G)) \approx 7 \cdot \frac{289}{71}+16 \cdot \frac{218}{71}+27 \cdot \frac{76}{71}+ 7 \cdot \frac{208}{71}+7 \cdot \frac{279}{71}+\frac{333.73}{71}+\frac{151.97}{71}+\frac{193.8}{71}+\frac{323.79}{71}+\frac{157.65}{71}+\frac{193.09}{71}+\frac{4736.38}{71} \approx \frac{17062.41}{71}$. Similarly, $|66-\frac{4752}{71}|=\frac{66}{71}, |\frac{129+\sqrt{33}}{2}-\frac{4752}{71}|\approx \frac{62.86}{142}, |\frac{129-\sqrt{33}}{2}-\frac{4752}{71}|=\frac{752.86}{142}$. Since $l_1 \approx 134.06, l_2 \approx 65.11$ and $l_3 \approx 55.82$, we have $|y_1-\frac{4752}{71}| \approx \frac{4766.26}{71}, |y_2 -\frac{4752}{71}| \approx \frac{129.19}{71}, |y_3 -\frac{4752}{71}| \approx \frac{788.78}{71}$ and hence $LE^{+}(\Gamma_{nc}(G)) \approx 27 \cdot \frac{76}{71}+ 18 \cdot \frac{66}{71}+7 \cdot \frac{279}{71}+8 \cdot \frac{62.86}{142}+8 \cdot \frac{752.86}{142}+ \frac{4766.26}{71}+\frac{129.19}{71}+\frac{788.78}{71} \approx \frac{28280.22}{142}$. 
	Further, $\CN-Spec(\Gamma_{nc}(G))=\left \lbrace (-68)^{27}, (-64)^{18}, (-63)^{7}, \left(\frac{-115-\sqrt{217}}{2} \right)^{8}, \left(\frac{-115+\sqrt{217}}{2} \right)^8, (m_1)^1, (m_2)^1, (m_3)^1 \right \rbrace$, where $m_1, m_2$ and $m_3$ are roots of the equation $x^3-4349x^2-311676x-1809504=0$. Since $m_1 \approx 4419.69, m_2 \approx -64.86$ and $m_3 \approx -6.37$ we have  $E_{CN}(\Gamma_{nc}(G)) \approx 8839.83$. Now, $|v(\Gamma_{nc}(G))|=71< 151.09=E(\Gamma_{nc}(G)), E(K_{71})=2(71-1)=140 < 151.09=E(\Gamma_{nc}(G))$ and $E_{CN}(K_{71})=2(71-1)(71-2)=9660>8839.83=E_{CN}(\Gamma_{nc}(G))$. Thus, \textit{$\Gamma_{nc}(G)$ is hyperenergetic but neither hypoenergetic nor CN-hyperenergetic}. Also, $LE(K_{71})=2(71-1)=140<\frac{17062.41}{71} \approx LE(\Gamma_{nc}(G))$ and $LE^{+}(K_{71})=2(71-1)=140<\frac{28280.22}{142} \approx LE^{+}(\Gamma_{nc}(G))$. Therefore, \textit{$\Gamma_{nc}(G)$ is L-hyperenergetic as well as Q-hyperenergetic}. Further, 
	\[
	E(\Gamma_{nc}(G)) \approx 151.09<\frac{28280.22}{142}=LE^{+}(\Gamma_{nc}(G))<\frac{17062.41}{71} \approx LE(\Gamma_{nc}(G)).
	\]
\end{proof}
\begin{theorem}
	Let $G$ be a finite non-abelian group such that $\Gamma_c(G)$ is triple-toroidal. Then 
	\begin{enumerate}
		\item $\Gamma_{nc}(G)$ is neither hypoenergetic, hyperenergetic nor CN-hyperenergetic.
		\item  $\Gamma_{nc}(G)$ is L-hyperenergetic as well as Q-hyperenergetic only when $G \cong GL(2, 3)$. 
		\item $E(\Gamma_{nc}(G)) \leq LE(\Gamma_{nc}(G)) \leq LE^{+}(\Gamma_{nc}(G))$.
	\end{enumerate}
\end{theorem}
\begin{proof}
	From Theorem \ref{graphs when triple toroidal}, we have that $\Gamma_c(G)$ is isomorphic to  $6K_2 \sqcup 3K_6 \sqcup 4K_4$ or  $3K_6$. 
	
	If $\Gamma_c(G) \cong 6K_2 \sqcup 3K_6 \sqcup 4K_4$, then $\Spec(\Gamma_{nc}(G))=\left \lbrace (0)^{33}, (-2)^5, (-6)^2, (-4)^3, (x_1)^1, \right .$
 $\left . (x_2)^1, (x_3)^1 \right \rbrace$, where $x_1, x_2$ and $x_3$ are roots of the equation $x^3-34x^2-312x-576=0$.
 Since  $x_1 \approx -5.08401 , x_2 \approx -2.71078, x_3 \approx 41.7948$, we have $E(\Gamma_{nc}(G)) \approx 10+12+12+5.08401+2.71078+41.7948 = 83.58959$. 
	We also have $\L-Spec(\Gamma_{nc}(G))=\left \lbrace (0)^{1}, (42)^{12}, (40)^{15}, (44)^{6}, (46)^{12} \right \rbrace$ and $\Q-Spec(\Gamma_{nc}(G))$ $=\left \lbrace (44)^6, (40)^{15}, (42)^{17}, (34)^2, (38)^{3}, (y_1)^1, (y_2)^1, (y_3)^1 \right \rbrace$, where $y_1, y_2$ and $y_3$ are roots of the equation $x^3-160x^2+7836x-121344=0$. Here, $\frac{2|e(\Gamma_{nc}(G))|}{|v(\Gamma_{nc}(G))|}=\frac{960}{23}$ so $|0-\frac{960}{23}|=\frac{960}{23}, |42-\frac{960}{23}|=\frac{6}{23}, |40-\frac{960}{23}|=\frac{40}{23}, |44-\frac{960}{23}|=\frac{52}{23}, |46-\frac{960}{23}|=\frac{98}{23}$. Therefore, $LE(\Gamma_{nc}(G))=\frac{960}{23}+12 \cdot \frac{6}{23}+15 \cdot \frac{40}{23}+ 6 \cdot \frac{52}{23}+12 \cdot \frac{98}{23}=\frac{3120}{23}$. Similarly, $|44-\frac{960}{23}|=\frac{52}{23}, |40-\frac{960}{23}|=\frac{40}{23}, |42-\frac{960}{23}|=\frac{6}{23}, |34-\frac{960}{23}|=\frac{178}{23}, |38-\frac{960}{23}|=\frac{86}{23}$. Since $y_1 \approx 35.7774, y_2 \approx 40.5202$ and $y_3 \approx 83.7024$, we have $|y_1-\frac{960}{23}| \approx 137.1198, |y_2 -\frac{960}{23}| \approx 28.0354, |y_3 -\frac{960}{23}| \approx 965.1552$ and hence $LE^{+}(\Gamma_{nc}(G)) \approx 6 \cdot \frac{52}{23}+15 \cdot \frac{40}{23}+17 \cdot \frac{6}{23}+2 \cdot \frac{178}{23}+ 3 \cdot \frac{86}{23}+137.1198+28.0354+965.1552=1201.0930$. 
	Further, $\CN-Spec(\Gamma_{nc}(G))=\left \lbrace (-44)^{6}, (-42)^{12}, (-40)^{20}, (-26)^{3}, (-4)^2, (z_1)^1, (z_2)^1, (z_3)^1 \right \rbrace$, where $z_1, z_2$ and $z_3$ are roots of the equation $x^3-1654x^2-86336x-921024=0$. Since $z_1 \approx 1704.96, z_2  \approx -35.9132$ and $z_3 \approx -15.042$ we have  $E_{CN}(\Gamma_{nc}(G)) \approx 3409.9152$. Now, $|v(\Gamma_{nc}(G))|=46< 83.58959=E(\Gamma_{nc}(G)), E(K_{46})=2(46-1)=90 > 83.58959=E(\Gamma_{nc}(G))$ and $E_{CN}(K_{46})=2(46-1)(46-2)=3960>3409.9152=E_{CN}(\Gamma_{nc}(G))$. Thus, \textit{$\Gamma_{nc}(G)$ is neither hypoenergetic, hyperenergetic nor CN-hyperenergetic}. Also, $LE(K_{46})=2(46-1)=90<\frac{3120}{23}=LE(\Gamma_{nc}(G))$ and $LE^{+}(K_{46})=2(46-1)=90<1201.0930=LE^{+}(\Gamma_{nc}(G))$. Therefore, \textit{$\Gamma_{nc}(G)$ is L-hyperenergetic as well as Q-hyperenergetic}. Further, 
	\[
	E(\Gamma_{nc}(G))=83.58959<\frac{3120}{23}=LE(\Gamma_{nc}(G))<1201.0930=LE^{+}(\Gamma_{nc}(G)).
	\]
	
	If $\Gamma_c(G) \cong 3K_6$, then 
	$\Spec(\Gamma_{nc}(G))=\left \lbrace (0)^{15}, (-6)^2, (12)^1 \right \rbrace$ and so $E(\Gamma_{nc}(G))=12+12=24$. 
	We also have $\L-Spec(\Gamma_{nc}(G))=\left \lbrace (0)^{1}, (12)^{15}, (18)^2 \right \rbrace$ and $\Q-Spec(\Gamma_{nc}(G))=\left \lbrace (6)^2, (12)^{15}, (24)^1 \right \rbrace$. Here, $\frac{2|e(\Gamma_{nc}(G))|}{|v(\Gamma_{nc}(G))|}=12$ so $|0-12|=12, |12-12|=0, |18-12|=6$. Therefore, $LE(\Gamma_{nc}(G))=12+0+2 \cdot 6=24$. Similarly, $|6-12|=6, |12-12|=0, |24-12|=12$ and hence $LE^{+}(\Gamma_{nc}(G))=2 \cdot 6 +0+ 12=24$. 
	Further, $\CN-Spec(\Gamma_{nc}(G))=\left \lbrace (132)^1, (24)^{2}, (-12)^{15} \right \rbrace$ and so $E_{CN}(\Gamma_{nc}(G))=360$. Now, $|v(\Gamma_{nc}(G))|=18< 24=E(\Gamma_{nc}(G)), E(K_{18})=2(18-1)=34 > 24=E(\Gamma_{nc}(G))$ and $E_{CN}(K_{18})=2(18-1)(18-2)=544>360=E_{CN}(\Gamma_{nc}(G))$. Thus, \textit{$\Gamma_{nc}(G)$ is neither hypoenergetic, hyperenergetic  nor CN-hyperenergetic}. Also, $LE(K_{18})=2(18-1)=34>24=LE(\Gamma_{nc}(G))$ and $LE^{+}(K_{18})=2(18-1)=34>24=LE^{+}(\Gamma_{nc}(G))$. Therefore, \textit{$\Gamma_{nc}(G)$ is neither L-hyperenergetic nor Q-hyperenergetic}. Further,
	\[
	E(\Gamma_{nc}(G))=24=LE(\Gamma_{nc}(G))=LE^{+}(\Gamma_{nc}(G)).
	\]
\end{proof}

{\bf Acknowledgement.} The first author is thankful to Council of Scientific and Industrial Research  for the fellowship (File No. 09/0796(16521)/2023-EMR-I).

\end{document}